\documentclass[12pt]{amsart}
\usepackage{amsmath}
\usepackage{amssymb}
\usepackage{amsfonts}
\usepackage{amsthm}
\usepackage{amscd}
\usepackage{bm}
\newtheorem{thrm}{Theorem}

\newtheorem{prop}{Proposition}
\newtheorem{thm}{Theorem}[section]
\newtheorem{lem}{Lemma}[section]
\newtheorem{pro}{Proposition}[section]
\newtheorem{cor}{Corollary}[section]

\newcommand{\Z}{\mathbf{Z}}

\begin{document}
\title[The kernel of the Magnus representation]
{The kernel of the Magnus representation of the automorphism group of a free group is not finitely generated}
\address{Department of Mathematics, Graduate School of Science, Kyoto University,
         Kitashirakawaoiwake cho, Sakyo-ku, Kyoto city 606-8502, Japan}
\email{takao@math.kyoto-u.ac.jp}
%
\maketitle
\begin{center}
{\sc Takao Satoh} \\

\vspace{0.5em}

{\footnotesize Department of Mathematics, Graduate School of Science, Kyoto University,
         Kitashirakawaoiwake-cho, Sakyo-ku, Kyoto city 606-8502, Japan}
\end{center}

\begin{abstract}
In this paper, we show that the abelianization of the kernel of the Magnus representation of the automorphism group of a free group
is not finitely generated.
\end{abstract}
\section{Introduction}\label{S-Int}

Let $F_n$ be a free group of rank $n \geq 2$, and $\mathrm{Aut}\,F_n$ the automorphism group of $F_n$.
Let denote $\rho : \mathrm{Aut}\,F_n \rightarrow \mathrm{Aut}\,H$
the natural homomorphism induced from the abelianization $F_n \rightarrow H$.
The kernel of $\rho$ is called the IA-automorphism group of $F_n$, denoted by $\mathrm{IA}_n$.
The subgroup $\mathrm{IA}_n$ reflects much richness and complexity of the structure of $\mathrm{Aut}\,F_n$,
and plays important roles on various studies of $\mathrm{Aut}\,F_n$.

\vspace{0.5em}

Although the study of the IA-automorphism group has a long history
since its finitely many generators were obtained by Magnus \cite{Mag} in 1935,
a presentation for $\mathrm{IA}_n$ is still unknown for $n \geq 3$.
Nielsen \cite{Ni0} showed that $\mathrm{IA}_2$ coincides with the inner automorphism group,
hence, is a free group of rank $2$.
For $n \geq 3$, however, $\mathrm{IA}_n$ is much larger than the inner automorphism group $\mathrm{Inn}\,F_n$.
Krsti\'{c} and McCool \cite{Krs} showed that $\mathrm{IA}_3$ is not finitely presentable.
For $n \geq 4$, it is not known whether $\mathrm{IA}_n$ is finitely presentable or not.

\vspace{0.5em}

In general, one of the most standard ways to study a group is to consider its representations.
If a representation of the group is given, it is important to determine whether it is faithful or not.
Furthermore, if not, it is also worth studying how far it is from faithful, namely, to determine its kernel.
In this paper, we consider these matters for the Magnus representation
\[ r_M : \mathrm{IA}_n \rightarrow \mathrm{GL}(n,\Z[H]) \]
of $\mathrm{IA}_n$. (See Subsection {\rmfamily \ref{Ss-Mag}}.)

\vspace{0.5em}

Classically, the Magnus representation was used to study certain subgroups of $\mathrm{IA}_n$.
One of the most famous subgroup is the pure braid group $P_n$. The restriction of $r_M$ to $P_n$ is called
the Gassner representation. (For a basic materials for the braid group and the pure braid group, see \cite{Bir} for example.)
It is known due to Magnus and Peluso \cite{MgP} that $r_M|_{P_n}$ is faithful for $n=3$.
Although the faithfulness of $r_M|_{P_n}$ has been studied for a long time by many authors,
it seems, however, still open problem to determine it for $n \geq 4$.

\vspace{0.5em}

Another important subgroup is the Torelli $\mathcal{I}_{g,1}$ subgroup of the mapping class group of a surface.
By a classical work of Dehn and Nielsen, the Torelli group $\mathcal{I}_{g,1}$ can be considered as a subgroup of $\mathrm{IA}_{2g}$.
Suzuki \cite{Su1} showed that the restriction $r_M|_{\mathcal{I}_{g,1}}$ to $\mathcal{I}_{g,1}$ is not faithful for $g \geq 2$.
Furthermore, by a recent remarkable work of Church and Farb \cite{ChF}, it is known that the abelianization of the kernel of
$r_M|_{\mathcal{I}_{g,1}}$ is not finitely generated for $g \geq 2$.

\vspace{0.5em}

From the facts as mentioned above, it is immediately seen that $r_M$ itself is not faithful. There are, however,
few results for the abelianization of the kernel $\mathcal{K}_n$ of $r_M$.
Let $F_n^M$ be the quotient group of $F_n$ by the second derived subgroup $[[F_n, F_n], [F_n, F_n]]$ of $F_n$.
The group $F_n^M$ is called the free metabelian group of rank $n$.
The metabelianization $F_n \rightarrow F_n^M$ naturally induces a homomorphism $\nu_n : \mathrm{Aut}\,F_n \rightarrow \mathrm{Aut}\,F_n^M$.
Its image of $\mathrm{IA}_n$ is contained in $\mathrm{IA}_n^M$, the IA-automorphism group of $F_n^M$.
Then it is known due to Bachmuth \cite{Ba1} that the Magnus representation $r_M$ factors through $\mathrm{IA}_n^M$, and that
the induced homomorphism $r_M': \mathrm{IA}_n^M \rightarrow \mathrm{GL}(n,\Z[H])$ is faithful.
Hence the metabelianization of $F_n$ induces the injectivization of the Magnus representation $r_M$.
Therefore we see that $\mathcal{K}_n$ coincides with the kernel of $\nu_n$.
From this point of view, in our previous paper \cite{S05}, we showed that the abelianization $\mathcal{K}_n^{\mathrm{ab}}$ of $\mathcal{K}_n$
contains a certain free abelian group of finite rank by using the Johnson homomorphism of $\mathrm{Aut}\,F_n^M$.

\vspace{0.5em}

In this paper, we show that $r_M$ is far from faithful. That is,
\begin{thrm}($=$ Theorem {\rmfamily \ref{T-main}} and Proposition {\rmfamily \ref{P-main}}.)
For $n \geq 2$, $\mathcal{K}_n^{\mathrm{ab}}$ is not finitely generated.
\end{thrm}

Recently, we have heard from Thomas Church, one of the authors of \cite{ChF}, that their method to show the infinite generation of
the abelianization of the kernel of $r_M|_{\mathcal{I}_{g,1}}$ can be applied to show that of $r_M$.
Our method, however, is different from theirs. In this paper, in order to prove the theorem for $n \geq 3$,
we consider some finitely generated normal subgroups $W_{n,d}$ of $F_n$, and
embeddings of $\mathcal{K}_n$ into the IA-automorphism group of $W_{n,d}$.
Then, taking advantage of the first Johnson homomorphisms of $\mathrm{Aut}\,W_{n,d}$, we detect infinitely many linearly independent elements in
$\mathcal{K}_n^{\mathrm{ab}}$. On the other hand, for $n=2$, we show
\begin{prop}
$\mathcal{K}_2 = [[F_2, F_2], [F_2, F_2]]$
\end{prop}
by using a usual and classical argument in combinatorial group theory.

\vspace{0.5em}

This paper consists of seven sections.
In Section {\rmfamily \ref{S-PRE}}, we recall the IA-automorphism group and the Magnus representation of the automorphism group of a free group.
In Section {\rmfamily \ref{S-Emb}}, we consider some embeddings of the kernel of the Magnus representation into the IA-automorphism group of $W_{n,d}$.
Then, in Section {\rmfamily \ref{S-Low}}, we prove the main theorem for $n \geq 3$.
In Section {\rmfamily \ref{S-det}}, we discuss how much of $\mathcal{K}_n^{\mathrm{ab}}$ can be detected in $\mathrm{IA}(W_{n,d})^{\mathrm{ab}}$.
In particular, we show that there exsists a non-trivial element in $\mathcal{K}_n^{\mathrm{ab}}$ which can not be detected by each of $\mathrm{IA}(W_{n,d})^{\mathrm{ab}}$
for $n \geq 4$.
In Section {\rmfamily \ref{S-n2}}, we consider the case where $n=2$.

\tableofcontents

\section{Preliminaries}\label{S-PRE}

In this section, after fixing notation and conventions,
we recall the IA-automorphism group and the Magnus representation of $\mathrm{Aut}\,F_n$.

\subsection{Notation and conventions}\label{Ss-NOT}
\hspace*{\fill}\ 

\vspace{0.5em}

Throughout the paper, we use the following notation and conventions. Let $G$ be a group and $N$ a normal subgroup of $G$.
\begin{itemize}
\item The abelianization of $G$ is denoted by $G^{\mathrm{ab}}$.
\item For any group $G$, we denote by $\Z[G]$ the integral group ring of $G$.
\item The group $\mathrm{Aut}\,G$ of $G$ acts on $G$ from the right.
      For any $\sigma \in \mathrm{Aut}\,G$ and $x \in G$, the action of $\sigma$ on $x$ is denoted by $x^{\sigma}$.
\item For an element $g \in G$, we also denote the coset class of $g$ by $g \in G/N$ if there is no confusion.
\item For elements $x$ and $y$ of $G$, the commutator bracket $[x,y]$ of $x$ and $y$
      is defined to be $[x,y]:=xyx^{-1}y^{-1}$.
\end{itemize}

\subsection{IA-automorphism group}\label{Ss-IA}
\hspace*{\fill}\ 

\vspace{0.5em}

In this paper, we fix a basis $x_1, \ldots , x_n$ of $F_n$.
Let $H:=F_n^{\mathrm{ab}}$ be the abelianization of $F_n$ and $\rho : \mathrm{Aut}\,F_n \rightarrow \mathrm{Aut}\,H$ the natural homomorphism induced from
the abelianization of $F_n$. In the following, we identify $\mathrm{Aut}\,H$ with the general linear group $\mathrm{GL}(n,\Z)$ by
fixing the basis of $H$ induced from the basis $x_1, \ldots , x_n$ of $F_n$.
The kernel $\mathrm{IA}_n$ of $\rho$ is called the IA-automorphism group of $F_n$.
It is clear that the inner automorphism group $\mathrm{Inn}\,F_n$
of $F_n$ is contained in $\mathrm{IA}_n$. In general, however, $\mathrm{IA}_n$ for $n \geq 3$ is much larger than $\mathrm{Inn}\,F_n$.
In fact, Magnus \cite{Mag} showed that for any $n \geq 3$, $\mathrm{IA}_n$ is finitely generated by automorphisms
\[ K_{ij} : x_t \mapsto \begin{cases}
               {x_j}^{-1} x_i x_j, & t=i, \\
               x_t,                & t \neq i
              \end{cases}\]
for distinct $i$, $j \in \{ 1, 2, \ldots , n \}$ and
\[  K_{ijl} : x_t \mapsto \begin{cases}
               x_i [x_j, x_l], & t=i, \\
               x_t,            & t \neq i
              \end{cases}\] 
for distinct $i$, $j$, $l \in \{ 1, 2, \ldots , n \}$ such that $j<l$.
In this paper, for the convenience, we also consider an automorphism $K_{ijl}$ for $j \geq l$ defined as above.
We remark that $K_{ijj}=1$ and $K_{ilj}=K_{ijl}^{-1}$.

\vspace{0.5em}

Recently, Cohen-Pakianathan \cite{Co1, Co2}CFarb \cite{Far} and Kawazumi \cite{Kaw} independently showed
\begin{equation}\label{CPFK}
\mathrm{IA}_n^{\mathrm{ab}} \cong H^* \otimes_{\Z} \Lambda^2 H
\end{equation}
as a $\mathrm{GL}(n,\Z)$-module where $H^*:= \mathrm{Hom}_{\Z}(H,\Z)$ is the $\Z$-linear dual group of $H$.

\vspace{0.5em}

In Section {\rmfamily \ref{S-Low}}, we use the following lemmas.
\begin{lem}\label{L-1}
For $n \geq 4$, and distinct integers $1 \leq i, j, l, m \leq n$, let $\omega \in \mathrm{IA}_n$ be an automorphism defined by
\[ x_t \mapsto \begin{cases}
                  x_i x_j x_l^{-1} x_m x_j^{-1} x_l x_m^{-1}, \hspace{1em} & t=i, \\
                  x_t, & t \neq i.
               \end{cases}\]
Then, we have
\[ \omega = K_{il} K_{iml} K_{ijm} K_{ilj} K_{il}^{-1} \in \mathrm{IA}_n, \]
and in particular,
\[ \omega = K_{iml} + K_{ijm} + K_{ilj} \in \mathrm{IA}_n^{\mathrm{ab}}. \]
\end{lem}
\textit{Proof.}
Since $K_{il} K_{iml} K_{ijm} K_{ilj} K_{il}^{-1}$ fix $x_t$ for $t \neq i$, it suffices to consider the image of $x_i$. Then we see
\[\begin{split}
   x_i & \xrightarrow{K_{il}} x_l^{-1} x_i x_l \xrightarrow{K_{iml}} x_l^{-1} x_i [x_m,x_l] x_l = x_l^{-1} x_i x_m x_l x_m^{-1}
       \xrightarrow{K_{ijm}} x_l^{-1} x_i [x_j, x_m] x_m x_l x_m^{-1} \\
       & = x_l^{-1} x_i x_j x_m x_j^{-1} x_l x_m^{-1} \xrightarrow{K_{ilj}} x_l^{-1} x_i [x_l, x_j] x_j x_m x_j^{-1} x_l x_m^{-1} \\
       & = x_l^{-1} x_i x_l x_j x_l^{-1} x_m x_j^{-1} x_l x_m^{-1} \xrightarrow{K_{il}^{-1}} x_i x_j x_l^{-1} x_m x_j^{-1} x_l x_m^{-1}.
  \end{split}\]
This completes the proof of Lemma {\rmfamily \ref{L-1}}. 

\vspace{0.5em}

Similarly, we obtain the following lemmas. Since the proofs are done by a straightforward calculation as above, we leave them to the reader
as exercises.

\begin{lem}\label{L-2}
For $n \geq 5$, and distinct integers $1 \leq i, j, l, m, p \leq n$, let $\omega \in \mathrm{IA}_n$ be an automorphism defined by
\[ x_t \mapsto \begin{cases}
                  x_i x_j x_l^{-1} x_p x_m x_p^{-1} x_j^{-1} x_l x_p x_m^{-1} x_p^{-1}, \hspace{1em} & t=i, \\
                  x_t, & t \neq i.
               \end{cases}\]
Then, we have
\[ \omega = K_{mp} K_{il} K_{iml} K_{ijm} K_{ilj} K_{il}^{-1} K_{mp}^{-1} \in \mathrm{IA}_n, \]
and in particular,
\[ \omega = K_{iml} + K_{ijm} + K_{ilj} \in \mathrm{IA}_n^{\mathrm{ab}}. \]
\end{lem}

\begin{lem}\label{L-4}
For $n \geq 5$, and distinct integers $1 \leq i, j, l, m, p \leq n$, let $\omega \in \mathrm{IA}_n$ be an automorphism defined by
\[ x_t \mapsto \begin{cases}
                  x_i x_p x_j x_p^{-1} x_l^{-1} x_p x_m x_p^{-1} x_p x_j^{-1} x_p^{-1}  x_l x_p x_m^{-1} x_p^{-1}, \hspace{1em} & t=i, \\
                  x_t, & t \neq i.
               \end{cases}\]
Then, we have
\[ \omega = K_{jp} K_{mp} K_{il} K_{iml} K_{ijm} K_{ilj} K_{il}^{-1} K_{mp}^{-1} K_{jp}^{-1} \in \mathrm{IA}_n, \]
and in particular,
\[ \omega = K_{iml} + K_{ijm} + K_{ilj} \in \mathrm{IA}_n^{\mathrm{ab}}. \]
\end{lem}

\subsection{Johnson homomorphisms}\label{Ss-Joh}
\hspace*{\fill}\ 

\vspace{0.5em}

Here we recall the definition of the Johnson homomorphisms of $\mathrm{Aut}\,F_n$.
(For a basic materials for the Johnson homomorphisms, see \cite{Mo1}, Kawazumi \cite{Kaw} and \cite{S03} for example.)

\vspace{0.5em}

Let $\Gamma_n(1) \supset \Gamma_n(2) \supset \cdots$ be the lower central series of a free group $F_n$ defined by the rule
\[ \Gamma_n(1):= F_n, \hspace{1em} \Gamma_n(k) := [\Gamma_n(k-1),F_n], \hspace{1em} k \geq 2. \]
We denote by $\mathcal{L}_n(k) := \Gamma_n(k)/\Gamma_n(k+1)$ the graded quotient of the lower central series of $F_n$.
For each $k \geq 1$, the action of $\mathrm{Aut}\,F_n$ on the nilpotent quotient group $F_n/\Gamma_n(k+1)$ of $F_n$ induces a homomorphism
$\mathrm{Aut}\,F_n \rightarrow \mathrm{Aut}(F_n/\Gamma_n(k+1))$.
We denote its kernel by $\mathcal{A}_n(k)$. Then the groups $\mathcal{A}_n(k)$ define a descending central filtration
\[ \mathrm{IA}_n = \mathcal{A}_n(1) \supset \mathcal{A}_n(2) \supset \cdots \]
of $\mathrm{IA}_n$. This filtration is called the Johnson filtration of $\mathrm{Aut}\,F_n$.

\vspace{0.5em}

The $k$-th Johnson homomorphism
\[ \tau_k : \mathcal{A}_n(k)/\mathcal{A}_n(k+1) \rightarrow \mathrm{Hom}_{\Z}(H, \mathcal{L}_n(k+1))
       = H^* \otimes_{\Z} \mathcal{L}_n(k+1) \]
of $\mathrm{Aut}\,F_n$ is defined by the rule
\[ \sigma \hspace{0.3em} \mapsto \hspace{0.3em} (x \mapsto x^{-1} x^{\sigma}), \hspace{1em} x \in H. \]
It is known that each of $\tau_k$ is a $\mathrm{GL}(n,\Z)$-equivariant injective homomorphism.
Since the target of the Johnson homomorphism $\tau_k$ is a free abelian group of finite rank, we can detect non-trivial elements
in the abelianization $\mathcal{A}_n(k)^{\mathrm{ab}}$ of $\mathcal{A}_n(k)$ by $\tau_k$.

\vspace{0.5em}

In particular, we remark that the isomorphism (\ref{CPFK}) is induced from the first Johnson homomorphism
\[ \tau_1 : \mathrm{IA}_n \rightarrow H^* \otimes_{\Z} \Lambda^2 H, \]
and that (the coset classes of) the Magnus generators $K_{ij}$s and $K_{ijk}$s form a basis of $\mathrm{IA}_n^{\mathrm{ab}}$
as a free abelian group.

\subsection{Magnus representation}\label{Ss-Mag}
\hspace*{\fill}\ 

\vspace{0.5em}

In this subsection we recall the Magnus representation of $\mathrm{Aut}\,F_n$. (For details, see \cite{Bir}.)
For each $1 \leq i \leq n$, let
\[ \frac{\partial}{\partial x_i} : \Z[F_n] \rightarrow \Z[F_n] \]
be the Fox derivation defined by
\[ \frac{\partial}{\partial x_i}(w) = 
                            \sum_{j=1}^{r} \epsilon_j \delta_{\mu_j,i} x_{\mu_1}^{\epsilon_1} \cdots
      x_{\mu_j}^{\frac{1}{2}(\epsilon_j-1)} \in \Z[F_n] \]
for any reduced word $w=x_{\mu_1}^{\epsilon_1} \cdots x_{\mu_r}^{\epsilon_r} \in F_n$, $\epsilon_j=\pm1$.
Let $\mathfrak{a} : F_n \rightarrow H$ be the abelianization of $F_n$. We also denote by $\mathfrak{a}$ the ring homomorphism
$\Z[F_n] \rightarrow \Z[H]$ induced from $\mathfrak{a}$.
For any matrix $A=(a_{ij}) \in \mathrm{GL}(n,\Z[F_n])$, set $A^{\mathfrak{a}} = (a_{ij}^{\mathfrak{a}}) \in \mathrm{GL}(n,\Z[H])$.
The Magnus representation $r_M : \mathrm{Aut}\,F_n \rightarrow \mathrm{GL}(n,\Z[H])$ of $\mathrm{Aut}\,F_n$
is defined by
\[ \sigma \mapsto \biggl{(} \frac{\partial {x_i}^{\sigma} }{\partial x_j} {\biggl{)}}^{\mathfrak{a}} \]
for any $\sigma \in \mathrm{Aut}\,F_n$.
This map is not a homomorphism but a crossed homomorphism. Namely,
\[ r_M(\sigma \tau) = r_M(\sigma)^{\tau^*} \cdot r_M(\tau) \]
where $r_M(\sigma)^{\tau^*}$ denotes the matrix obtained from $r_M(\sigma)$
by applying the automorphism
$\tau^* : \Z[H] \rightarrow \Z[H]$ induced from $\rho(\tau) \in \mathrm{Aut}(H)$ on each entry.
Hence by restricting $r_M$ to $\mathrm{IA}_n$, we obtain a homomorphism
$\mathrm{IA}_n \rightarrow \mathrm{GL}(n,\Z[H])$, also denote it by $r_M$.

\vspace{0.5em}

Let $\mathcal{K}_n$ be the kernel of the homomorphism $r_M : \mathrm{IA}_n \rightarrow \mathrm{GL}(n,\Z[H])$.
The main purpose of the paper is to show that the abelianization $\mathcal{K}_n^{\mathrm{ab}}$ of $\mathcal{K}_n$ is not finitely generated.
More precisely, we prove that it contains a free abelian group of infinite rank.

\vspace{0.5em}

Here we consider the group $\mathcal{K}_n$ from the view point of the metabelianization of $F_n$.
Let $F_n^M$ be the quotient group of $F_n$ by the second derived subgroup $[[F_n, F_n], [F_n, F_n]]$ of $F_n$.
The group $F_n^M$ is called the free metabelian group of rank $n$. The abelianization of $F_n^M$ is canonically isomorphic to $H=F_n^{\mathrm{ab}}$.
Let $\mu_n : \mathrm{Aut}\,F_n \rightarrow \mathrm{Aut}\,F_n^M$ be the induced homomorphism from the action of $\mathrm{Aut}\,F_n$ on $F_n^M$.
It is known that $\mu_n$ is surjective for $n \geq 2$ except for $n=3$. (See \cite{Ba1} for $n=2$, and see \cite{BM2} for $n \geq 4$.)
Let $\mathrm{IA}_n^M$ be the IA-automorphism group of $F_n^M$. Namely, the group $\mathrm{IA}_n^M$ consists of automorphisms of $F_n^M$ which act on
the abelianization of $F_n^M$ trivially. Then the restriction of $\mu_n$ induces a homomorphism $\mathrm{IA}_n \rightarrow \mathrm{IA}_n^M$,
also denoted by $\mu_n$.

\vspace{0.5em}

Now, in \cite{Ba1}, Bachmuth constructed a faithful representation
\[ r_M' : \mathrm{IA}_n^M \rightarrow \mathrm{GL}(n,\Z[H]) \]
using the Magnus representation of $F_n^M$.
Then we can easily see that $r_M = r_M' \circ \mu_n$ by observing the image of the Magnus generators of $\mathrm{IA}_n$.
(For details, see \cite{Ba1}.) Therefore, the faithfulness of $r_M'$ induces $\mathcal{K}_n = \mathrm{Ker}(\mu_n)$. In particular, we have
\begin{equation}
 \mathcal{K}_n = \{ \sigma \in \mathrm{Aut}\,F_n \,|\, x^{-1} x^{\sigma} \in [[F_n, F_n], [F_n, F_n]], \,\,\, x \in F_n \}
   \subset \mathrm{IA}_n. \label{eq-ker}
\end{equation}

\section{Embeddings of the IA-automorphism group}\label{S-Emb}

In this section, we consider certain finitely generated subgroups of the free group $F_n$, and its automorphism group.

\subsection{Subgroups $W_{n,d}$ of $F_n$}\label{Ss-Sub}
\hspace*{\fill}\ 

\vspace{0.5em}

For any integer $d \geq 2$, let $C_d$ be the cyclic group of order $d$ generated by $s$. Let $\iota_{n,d}: F_n \rightarrow C_d$
be a homomorphism defined by
\[ x_i^{\iota_{n,d}} := \begin{cases}
                     s, \hspace{1em} & i=1, \\
                     1, & i \neq 1.
                  \end{cases} \]
We denote by $W_{n,d}$ the kernel of $\iota_{n,d}$. Then, $W_{n,d}$ is also a free group, and its rank is given by $d(n-1)+1$
since the index of $W_{n,d}$ in $F_n$ is $d$. Furthermore, applying the Reidemeister-Schreier method to the following data:
\[\begin{split}
   X & := \{ x_1, x_2, \ldots, x_n \}; \,\,\, \text{the set of generators of $F_n$}, \\
   T & := \{ 1, x_1, \ldots, x_1^{d-1} \}; \,\,\, \text{a Schreier-transversal of $W_{n,d}$ of $F_n$}, \\
  \end{split}\]
we see that a set
\begin{equation}\label{eq-1}
 \{ (t,x) \in W_{n,d} \,|\, t \in T, \,\,\, x \in X, \,\,\, (t,x) \neq 1 \}
\end{equation}
form a basis of $W_{n,d}$, where $(t,x)=tx(\overline{tx})^{-1}$ and a map $\bar{} : F_n \rightarrow T$ is defined by the condition that
\[ W_{n,d} \,\, y = W_{n,d} \,\, \overline{y} \subset F_n \]
for any $y \in F_n$. (For details for the Reidemeister-Schreier method, see \cite{Lyn} for example.)

\vspace{0.5em}

Then the set (\ref{eq-1}) is written down as
\begin{equation}\label{eq-b}
\begin{split}
   & x_1^d, \\
   & x_2, \, x_1 x_2 x_1^{-1}, \, \ldots, \, x_1^{d-1} x_2 x_1^{-(d-1)}, \\
   & \hspace{6em} \vdots \\
   & x_n, \, x_1 x_n x_1^{-1}, \, \ldots, \, x_1^{d-1} x_n x_1^{-(d-1)}. \\
  \end{split}\end{equation}
In the following, we fix this set as a basis of $W_{n,d}$.

\vspace{0.5em}

\subsection{Action of $\mathrm{IA}_n$ on $W_{n,d}$}\label{Ss-Act}
\hspace*{\fill}\ 

\vspace{0.5em}

Since each of $W_{n,d}$ contains the commutator subgroup $[F_n,F_n]$ of $F_n$, the IA-automorphism group
\[ \mathrm{IA}_n = \{ \sigma \in \mathrm{Aut}\,F_n \,|\, x^{-1} x^{\sigma} \in [F_n,F_n], \,\,\, x \in F_n \} \]
naturally acts on them. Let $\rho = \rho_{n,d} : \mathrm{IA}_n \rightarrow \mathrm{Aut}\,W_{n,d}$ be the homomorphism induced from the action.
Then we have

\begin{lem}
For any $n \geq 2$ and $d \geq 2$, the homomorphism $\rho$ is injective.
\end{lem}
\textit{Proof.}
For an element $\sigma \in \mathrm{IA}_n$, assume $\rho(\sigma)=1$. Since $x_i^{\sigma} = x_i$ for $2 \leq i \leq n$, in order to see $\sigma=1$,
it suffices to show $x_1^{\sigma}=x_1$.

\vspace{0.5em}

Set $x_1^{\sigma}=x_1 v$ for some irreducible word $v$. Then we have
\begin{equation}
 x_1^d = (x_1^d)^{\sigma} = (x_1^{\sigma})^d = x_1 v x_1 v \cdots x_1 v, \label{eq-2}
\end{equation}
In the right hand side of the equation above, no cancellation of words happen. In fact, if $v$ is of type $x_1^{-1} w$ for some irreducible
word $w$, then we have $x_1^d= w^d$. By a usual argument in the Combinatorial group theory, there exists some $z \in F_n$ such that
$x_1 =z^a$ and $w=z^b$ for some $a, b \in \Z$. (See proposition 2.17. in \cite{Lyn}.) From the former equation, we have $z=x_1$ and $a=1$.
Hence
\[ x_1^d \, w^{-d} = x_1^{d - bd} = 1. \]
Since $F_n$ is torsion free, $b=1$. This shows that $w=x_1$. This is a contradiction to the irreducibility of $v=x_1^{-1} w$.
On the other hand, If $v$ is of type $w x_1^{-1}$ for some irreducible word $w$, then we have
\[ x_1^d = x_1^{-1} \, x_1^d \, x_1 = x_1^{-1} \, (x_1 w^d x_1^{-1}) x_1 = w^d. \]
By an argument similar to the above, we see that this is also contradiction.

\vspace{0.5em}

Thus, observing the word length of both side of (\ref{eq-2}), we obtain that $v=1$. This shows that $\sigma=1$. 

\vspace{0.5em}

By this lemma, we can consider $\mathrm{IA}_n$ as a subgroup of each of $\mathrm{Aut}\,W_{n,d}$. In the following, we identify $\mathrm{IA}_n$
with the image of it by $\rho$.

\vspace{1em}

Now, let $U_{n,d}$ be the abelianization of $W_{n,d}$. We denote by $\mathrm{IA}(W_{n,d})$ the IA-automorphism group of $W_{n,d}$.
Namely, the group $\mathrm{IA}(W_{n,d})$ consists of automorphisms of $W_{n,d}$ which act on $U_{n,d}$ trivially.
Let us define an linear order among (\ref{eq-b}) by
\[ x_1^d < x_1 x_2 x_1^{-1} < \ldots < x_1^{d-1} x_2 x_1^{-(d-1)} < \ldots < x_1 x_n x_1^{-1} < \ldots < x_1^{d-1} x_n x_1^{-(d-1)}. \]
Then, by a result of Magnus \cite{Mag} as mentioned in Subsection {\rmfamily \ref{Ss-IA}}, the IA-automorphism group $\mathrm{IA}(W_{n,d})$ is
finitely generated by automorphisms
\begin{equation}
\begin{split}
 & K_{\alpha, \beta}: \alpha \mapsto \alpha^{-1} \beta \alpha, \hspace{1em} \alpha < \beta, \\
 & K_{\alpha, \beta, \gamma}: \alpha \mapsto \alpha [\beta, \gamma], \hspace{1em} \alpha \neq \beta < \gamma \neq \alpha
\end{split}
\end{equation}
where $\alpha$, $\beta$, $\gamma$ are elements of the basis (\ref{eq-b}).

\vspace{0.5em}

On the other hand, from (\ref{CPFK}), we have
\[ \mathrm{IA}(W_{n,d})^{\mathrm{ab}} \cong U_{n,d}^* \otimes_{\Z} \Lambda^2 U_{n,d} \]
where $U_{n,d}^*$ denotes the $\Z$ linear dual group $\mathrm{Hom}_{\Z}(U_{n,d},\Z)$ of $U_{n,d}$.
This isomorphism is induced from the first Johnson homomorphism
\[ \tau_1 : \mathrm{IA}(W_{n,d}) \rightarrow U_{n,d}^* \otimes_{\Z} \Lambda^2 U_{n,d}. \]
In particular, (the coset classes of) the elements $K_{\alpha, \beta}$s and $K_{\alpha, \beta, \gamma}$s form a basis of $\mathrm{IA}(W_{n,d})^{\mathrm{ab}}$
as a free abelian group.

\vspace{0.5em}

In general, the induced action of $\mathrm{IA}_n$ on $U_{n,d}$ by $\rho$ is not trivial. However, its restriction to the kernel
of the Magnus representation $\mathcal{K}_n$ is trivial by (\ref{eq-ker}).
Namely, $\rho$ induces an embedding
\[ \rho|_{\mathcal{K}_n} : \mathcal{K}_n \hookrightarrow \mathrm{IA}(W_{n,d}). \]
By this embedding, we consider $\mathcal{K}_n$ as a subgroup of $\mathrm{IA}(W_{n,d})$ for any $d \geq 2$.

\section{Lower bounds on $\mathcal{K}_n^{\mathrm{ab}}$}\label{S-Low}

In the following, we always assume $n \geq 3$ and $d \geq 3$.
For each $2 \leq m \leq d-1$, let $\sigma_m \in \mathrm{Aut}\,F_n$ be an automorphism defined by
\[ x_i \mapsto \begin{cases}
           x_2 [[x_1, x_3], [x_1^m, x_3]], \hspace{1em} & i=2, \\
           x_i, & i \neq 2.
           \end{cases} \]
Clearly, we see $\sigma_m \in \mathcal{K}_n$. We show that $\sigma_2, \ldots, \sigma_{d-1}$ are linearly independent in $\mathcal{K}_n^{\mathrm{ab}}$.
To do this, we consider the images of $\sigma_m$ by the natural homomorphism
\[ \pi_{n,d} : \mathcal{K}_n \hookrightarrow \mathrm{IA}(W_{n,d}) \rightarrow \mathrm{IA}(W_{n,d})^{\mathrm{ab}}. \]

\vspace{0.5em}

To begin with, in order to describe $\sigma_m$ in $\mathrm{IA}(W_{n,d})^{\mathrm{ab}}$ with $K_{\alpha, \beta}$ and $K_{\alpha, \beta, \gamma}$,
we consider the action of $\sigma_m$ on the basis (\ref{eq-b}) of $W_{n,d}$.
Observing
\[ x_2^{\sigma_m} = x_2 \cdot x_1x_3 x_1^{-1} \cdot x_3^{-1} \cdot x_1^m x_3 x_1^{-m} \cdot x_1 x_3^{-1} x_1^{-1} \cdot x_3
                    \cdot x_1^m x_3^{-1} x_1^{-m}, \]
we see
\[\begin{split}
   (x_1^d)^{\sigma_m} & = x_1^d, \\
   (x_1^k x_i x_1^{-k})^{\sigma_m} & = x_1^k x_i x_1^{-k}, \hspace{1em} i \neq 2, \\
   (x_1^k x_2 x_1^{-k})^{\sigma_m} & = x_1^k x_2 x_1^{-k} \cdot x_1^{k+1} x_3 x_1^{-(k+1)} \cdot x_1^k x_3^{-1} x_1^{-k}
    \cdot x_1^{m+k} x_3 x_1^{-(m+k)} \\
    & \hspace{3em} \cdot x_1^{k+1} x_3^{-1} x_1^{-(k+1)} \cdot x_1^k x_3 x_1^{-k} \cdot x_1^{m+k} x_3^{-1} x_1^{-(m+k)}.
 \end{split}\]
Hence, for any $0 \leq k \leq d-1$, if we denote by $\omega_{m,k}$ the automorphism of $W_{n,d}$ defined by
\begin{equation}\label{eq-i}
\begin{split}
   x_1^k x_2 x_1^{-k} & \mapsto x_1^k x_2 x_1^{-k} \cdot x_1^{k+1} x_3 x_1^{-(k+1)} \cdot x_1^k x_3^{-1} x_1^{-k}
    \cdot x_1^{m+k} x_3 x_1^{-(m+k)} \\
   & \hspace{3em} \cdot x_1^{k+1} x_3^{-1} x_1^{-(k+1)} \cdot x_1^k x_3 x_1^{-k} \cdot x_1^{m+k} x_3^{-1} x_1^{-(m+k)}, 
  \end{split}
\end{equation}
then
\[ \sigma_m = \omega_{m,0} \omega_{m,2} \cdots \omega_{m, d-1}. \]
Here we remark that for any $0 \leq k, l \leq d-1$, the automorphisms $\omega_{m,k}$ and $\omega_{m,l}$ are commutative in $\mathrm{Aut}\,W_{n,d}$.

\vspace{0.5em}

Let us describe $\omega_{m,k}$ in $\mathrm{IA}(W_{n,d})^{\mathrm{ab}}$ with $K_{\alpha, \beta}$ and $K_{\alpha, \beta, \gamma}$ for each $0 \leq k \leq d-1$.

\begin{flushleft}
{\bf Case 1.} The case of $0 \leq k \leq d-1-m$.

\vspace{0.5em}

In this case, we have $k+1 \leq d-1$ and $m+k \leq  d-1$. Thus, by (\ref{eq-i}) and Lemma {\rmfamily \ref{L-1}}, we see
\[\begin{split}
   \omega_{m,k} & = K_{x_1^k x_2 x_1^{-k}, \, x_1^k x_3 x_1^{-k}} K_{x_1^k x_2 x_1^{-k}, \, x_1^{m+k} x_3 x_1^{-(m+k)}, \, x_1^k x_3 x_1^{-k}} \\
                & \hspace{2em} \cdot K_{x_1^k x_2 x_1^{-k}, \, x_1^{k+1} x_3 x_1^{-(k+1)}, \, x_1^{m+k} x_3 x_1^{-(m+k)}} \\
                & \hspace{2em} \cdot K_{x_1^k x_2 x_1^{-k}, \, x_1^k x_3 x_1^{-k}, \, x_1^{k+1} x_3 x_1^{-(k+1)}}
                    K_{x_1^k x_2 x_1^{-k}, \, x_1^k x_3 x_1^{-k}}^{-1}
  \end{split}\]
in $\mathrm{IA}(W_{n,d})$, and hence,
\[\begin{split}
    \omega_{m,k} & = K_{x_1^k x_2 x_1^{-k}, \, x_1^{m+k} x_3 x_1^{-(m+k)}, \, x_1^k x_3 x_1^{-k}}
              + K_{x_1^k x_2 x_1^{-k}, \, x_1^{k+1} x_3 x_1^{-(k+1)}, \, x_1^{m+k} x_3 x_1^{-(m+k)}} \\
           & \hspace{2em} + K_{x_1^k x_2 x_1^{-k}, \, x_1^k x_3 x_1^{-k}, \, x_1^{k+1} x_3 x_1^{-(k+1)}}. 
  \end{split}\]
in $\mathrm{IA}(W_{n,d})^{\mathrm{ab}}$.
\end{flushleft}

\begin{flushleft}
{\bf Case 2.} The case of $d-m \leq k \leq d-2$.

\vspace{0.5em}

In this case, we have $k+1 \leq d-1 < m+k$. Hence (\ref{eq-i}) is written as
\[\begin{split}
   x_1^k x_2 x_1^{-k} & \mapsto x_1^k x_2 x_1^{-k} \cdot x_1^{k+1} x_3 x_1^{-(k+1)} \cdot x_1^k x_3^{-1} x_1^{-k}
    \cdot x_1^d \cdot x_1^{m+k-d} x_3 x_1^{-(m+k-d)} \cdot x_1^{-d} \\
   & \hspace{3em} \cdot x_1^{k+1} x_3^{-1} x_1^{-(k+1)} \cdot x_1^k x_3 x_1^{-k} \cdot x_1^d \cdot x_1^{m+k-d} x_3^{-1} x_1^{-(m+k-d)} \cdot x_1^{-d}.
  \end{split}\]
Then using Lemma {\rmfamily \ref{L-2}}, we obtain
\[\begin{split}
    \omega_{m,k} & = K_{x_1^k x_2 x_1^{-k}, \, x_1^{m+k-d} x_3 x_1^{-(m+k-d)}, \, x_1^k x_3 x_1^{-k}}
              + K_{x_1^k x_2 x_1^{-k}, \, x_1^{k+1} x_3 x_1^{-(k+1)}, \, x_1^{m+k-d} x_3 x_1^{-(m+k-d)}} \\
           & \hspace{2em} + K_{x_1^k x_2 x_1^{-k}, \, x_1^k x_3 x_1^{-k}, \, x_1^{k+1} x_3 x_1^{-(k+1)}}. 
  \end{split}\]
in $\mathrm{IA}(W_{n,d})^{\mathrm{ab}}$.
\end{flushleft}

\begin{flushleft}
{\bf Case 3.} The case of $k=d-1$.

\vspace{0.5em}

In this case, we have $k+1=d$ and $d \leq m+k$. Hence (\ref{eq-i}) is written as
\[\begin{split}
   x_1^{d-1} x_2 x_1^{-(d-1)} & \mapsto x_1^{d-1} x_2 x_1^{-(d-1)} \\
   & \hspace{3em} \cdot x_1^{d} \cdot x_3 \cdot x_1^{-d} \cdot x_1^{d-1} x_3^{-1} x_1^{-(d-1)}
    \cdot x_1^d \cdot x_1^{m-1} x_3 x_1^{-(m-1)} \cdot x_1^{-d} \\
   & \hspace{3em} \cdot x_1^d \cdot x_3^{-1} \cdot x_1^{-d} \cdot x_1^{d-1} x_3 x_1^{-(d-1)} \cdot x_1^d \cdot x_1^{m-1} x_3^{-1} x_1^{-(m-1)} \cdot x_1^{-d}.
  \end{split}\]
Then using Lemma {\rmfamily \ref{L-4}}, we obtain
\[\begin{split}
    \omega_{m,d-1} & = K_{x_1^{d-1} x_2 x_1^{-(d-1)}, \, x_1^{m-1} x_3 x_1^{-(m-1)}, \, x_1^{d-1} x_3 x_1^{-(d-1)}}
              + K_{x_1^{d-1} x_2 x_1^{-(d-1)}, \, x_3, \, x_1^{m-1} x_3 x_1^{-(m-1)}} \\
           & \hspace{2em} + K_{x_1^{d-1} x_2 x_1^{-(d-1)}, \, x_1^{d-1} x_3 x_1^{-(d-1)}, \, x_3}. 
  \end{split}\]
in $\mathrm{IA}(W_{n,d})^{\mathrm{ab}}$.
\end{flushleft}

\vspace{0.5em}

Therefore for $2 \leq m \leq d-1$, we obtain
\begin{small}
\begin{equation}\begin{split} \label{eq-mayumi}
   \sigma_m & = \sum_{0 \leq k \leq d-1-m} \{ K_{x_1^k x_2 x_1^{-k}, \, x_1^{m+k} x_3 x_1^{-(m+k)}, \, x_1^k x_3 x_1^{-k}} \\
            & \hspace{3em} + K_{x_1^k x_2 x_1^{-k}, \, x_1^{k+1} x_3 x_1^{-(k+1)}, \, x_1^{m+k} x_3 x_1^{-(m+k)}}
               + K_{x_1^k x_2 x_1^{-k}, \, x_1^k x_3 x_1^{-k}, \, x_1^{k+1} x_3 x_1^{-(k+1)}} \} \\
            & \hspace{1em} + \sum_{d-m \leq k \leq d-1} \{ K_{x_1^k x_2 x_1^{-k}, \, x_1^{m+k-d} x_3 x_1^{-(m+k-d)}, \, x_1^k x_3 x_1^{-k}} \\
            & \hspace{3em} + K_{x_1^k x_2 x_1^{-k}, \, x_1^{k+1} x_3 x_1^{-(k+1)}, \, x_1^{m+k-d} x_3 x_1^{-(m+k-d)}}
               + K_{x_1^k x_2 x_1^{-k}, \, x_1^k x_3 x_1^{-k}, \, x_1^{k+1} x_3 x_1^{-(k+1)}} \}  \\
  \end{split}\end{equation}
\end{small}
in $\mathrm{IA}(W_{n,d})^{\mathrm{ab}}$. Using this, we obtain
\begin{lem}\label{L-ml}
For any $n \geq 3$ and $d \geq 3$, the elements $\sigma_2, \ldots, \sigma_{d-1}$ are linearly independent in $\mathrm{IA}(W_{n,d})^{\mathrm{ab}}$.
\end{lem}
\textit{Proof.}
Suppose
\[ a_2 \sigma_2 + \cdots + a_d \sigma_d = 0 \in \mathrm{IA}(W_{n,d})^{\mathrm{ab}}. \]
Then, observing the coefficients of $K_{x_2, *, *}$ in the left hand side of (\ref{eq-mayumi}), we see
\[\begin{split}
   & \sum_{2 \leq m \leq d-1} a_m \{ K_{x_2, \, x_1^{m} x_3 x_1^{-m}, \, x_3} + K_{x_2, \, x_1 x_3 x_1^{-1}, \, x_1^{m} x_3 x_1^{-m}}
               + K_{x_2, \, x_3, \, x_1 x_3 x_1^{-1}} \} = 0.
   \end{split}\]
Hence, we have $a_2= \cdots = a_{d-1}=0$. 

\vspace{0.5em}

This Lemma induces our main theorem.
\begin{thm}\label{T-main}
For $n \geq 3$, $\mathcal{K}_n^{\mathrm{ab}}$ is not finitely generated.
\end{thm}
\textit{Proof.}
For any $d \geq 3$, the image of the induced map
\[ \mathcal{K}_n^{\mathrm{ab}} \rightarrow \mathrm{IA}(W_{n,d})^{\mathrm{ab}} \]
from the natural homomorphism $\pi_{n,d} : \mathcal{K}_n \rightarrow \mathrm{IA}(W_{n,d})^{\mathrm{ab}}$
contains a free abelian group of rank $d-2$. Since we can take $d \geq 3$ arbitrarily, we obtain the required result.
This completes the proof of Theorem {\rmfamily \ref{T-main}}. 

\vspace{0.5em}

As a Corollary, we have

\begin{cor}
For $n \geq 3$, $\mathcal{K}_n$ is not finitely generated.
\end{cor}

\section{On the detectivity of elements of $\mathcal{K}_n^{\mathrm{ab}}$ in $\mathrm{IA}(W_{n,d})^{\mathrm{ab}}$}\label{S-det}

In this section, we consider how much of $\mathcal{K}_n^{\mathrm{ab}}$ can be detected in $\mathrm{IA}(W_{n,d})^{\mathrm{ab}}$.
In the previous section, we have shown that $\mathrm{IA}(W_{n,d})^{\mathrm{ab}}$ detects at least $\Z^{\oplus d-2}$ in $\mathcal{K}_n^{\mathrm{ab}}$.
First, we give a more sharp result for the detectivity of linearly independent elements of $\mathcal{K}_n^{\mathrm{ab}}$
in $\mathrm{IA}(W_{n,d})^{\mathrm{ab}}$.
Second, we show that there is a non-trivial element in $\mathcal{K}_n^{\mathrm{ab}}$ which can not be detected in all $\mathrm{IA}(W_{n,d})^{\mathrm{ab}}$
for $n \geq 4$.

\vspace{0.5em}

For any $d \geq 3$ and $2 \leq m \leq d-1$, let $\sigma_m^{j,s}$ be an automorphism of $F_n$ defined by
\[ x_i \mapsto \begin{cases}
           x_j [[x_1, x_s], [x_1^m, x_s]], \hspace{1em} & i=j, \\
           x_i, & i \neq j
           \end{cases}\]
for $2 \leq j \leq n$, $2 \leq s \leq n$ and $s \neq j$. Then we see that $\sigma_m^{j,s} \in \mathcal{K}_n$ and $\sigma_m^{2,3} = \sigma_m$.

\begin{lem}\label{L-TY}
For any $n \geq 3$ and $d \geq 3$, the elements $\sigma_m^{j,s}$ are linearly independent in $\mathrm{IA}(W_{n,d})^{\mathrm{ab}}$.
\end{lem}
\textit{Proof.}
Suppose
\[ \sum_{m=2}^{d-1} \sum_{j=2}^n \sum_{2 \leq s \leq n}^{s \neq j} a_m^{j,s} \sigma_m^{j,s} =0
   \in \mathrm{IA}(W_{n,d})^{\mathrm{ab}}. \]
By the same argument as (\ref{eq-mayumi}), $\sigma_m^{j,s} \in \mathrm{IA}(W_{n,d})^{\mathrm{ab}}$ is a linear combination
of the elements type of $K_{x_1^k x_j x_1^{-k}, *, *}$. Hence, for each $2 \leq j \leq n$,
\[ \sum_{m=2}^{d-1} \sum_{2 \leq s \leq n}^{s \neq j} a_m^{j,s} \sigma_m^{j,s} =0. \]
Furthermore, for each $s \neq j$, since $\sigma_m^{j,s} \in \mathrm{IA}(W_{n,d})^{\mathrm{ab}}$ is a linear combination
of the elements type of $K_{x_1^k x_j x_1^{-k}, x_1^a x_s x_1^{-a}, x_1^b x_s x_1^{-b}}$, we see
\[ \sum_{m=2}^{d-1} a_m^{j,s} \sigma_m^{j,s} =0. \]
Fianlly, if we obverve the coefficients of $K_{x_j, *, *}$ in the equation above, we see
\[ \sum_{m=2}^{d-1} a_m^{j,s} \{ K_{x_j, \, x_1^{m} x_s x_1^{-m}, \, x_s} + K_{x_j, \, x_1 x_s x_1^{-1}, \, x_1^{m} x_s x_1^{-m}}
               + K_{x_j, \, x_s, \, x_1 x_s x_1^{-1}} \} = 0. \]
Therefore, we obtain $a_m^{j,s}=0$. This completes the proof of Lemma {\rmfamily \ref{L-TY}}. 

\vspace{0.5em}

As a corollary, we see that $\mathrm{IA}(W_{n,d})^{\mathrm{ab}}$ detects $\Z^{\oplus (d-2)(n-1)(n-2)}$ in $\mathcal{K}_n^{\mathrm{ab}}$.
We remark that in general, it seems difficult to determine the image of the homomorphism
$\mathcal{K}_n^{\mathrm{ab}} \rightarrow \mathrm{IA}(W_{n,d})^{\mathrm{ab}}$.

\vspace{1em}

Next, we construct a non-trivial element in $\mathcal{K}_n^{\mathrm{ab}}$ which can not be detected in all $\mathrm{IA}(W_{n,d})^{\mathrm{ab}}$.
In the following, we assume $n \geq 4$ until the end of this section.
For $1 \leq i \leq n$, we denote by $\iota_i \in \mathrm{IA}_n$ the inner automorphism of $F_n$ defined by
\[ x \mapsto x_i x x_i^{-1}, \hspace{1em} x \in F_n. \]
Set
\[ \sigma = [[\iota_2, \iota_3], [\iota_2, \iota_4]] \in \mathrm{IA}_n. \]
Then we have
\[ x^{\sigma} = [[x_2, x_3], [x_2, x_4]] \, x \, [[x_2, x_4], [x_2, x_3]], \]
and $\sigma \in \mathcal{K}_n$.

\begin{lem}\label{L-kouko}
$\sigma \neq 0 \in \mathcal{K}_n^{\mathrm{ab}}$.
\end{lem}
\textit{Proof.}
Suppose $\sigma \in [\mathcal{K}_n, \mathcal{K}_n]$. By (\ref{eq-ker}), we see that $\mathcal{K}_n \subset \mathcal{A}_n(3)$.
Since the Johnson filtration is central, we have $\sigma \in \mathcal{A}_n(6)$.

\vspace{0.5em}

On the other hand, observing the image of $\sigma$ by $\tau_4$, we have
\[ \tau_4 (\sigma) = \sum_{i=1}^n x_i^* \otimes [x_i^{-1}, [[x_2, x_3], [x_2, x_4]]] \in H^* \otimes_{\Z} \mathcal{L}_n(5) \]
where $x_i^*$s are the dual basis of $x_i$s. In particular, for each $i \geq 2$ we have
\[\begin{split}
     [x_i^{-1}, & [[x_2, x_3], [x_2, x_4]]] = [[[x_2, x_3], [x_2, x_4]], x_i], \\
                & = -[[[x_2, x_4], x_i], [x_2, x_3]] + [[[x_2, x_3], x_i], [x_2, x_4]], \hspace{1em} (\text{Jacobi identity}), \\
                & = [[[x_4, x_2], x_i], [x_2, x_3]] - [[[x_3, x_2], x_i], [x_2, x_4]], \\
                & \neq 0 \in \mathcal{L}_n(5)
  \end{split}\]
since both of
\[ [[[x_4, x_2], x_i], [x_2, x_3]] \,\,\, \text{and} \,\,\, [[[x_3, x_2], x_i], [x_2, x_4]] \]
are elements of the Hall basis of $\mathcal{L}_n(5)$ with the standard ordering $x_1 < x_2 < \cdots < x_n$.
(For a basic materials for the Hall basis, see \cite{Reu} for example.) This is a contradiction. This completes the proof of
Lemma {\rmfamily \ref{L-kouko}}. 

\vspace{0.5em}

It is clear that $\sigma =0 \in \mathrm{IA}(W_{n,d})^{\mathrm{ab}}$ for any $d \geq 3$ since both of
$[\iota_2, \iota_3]$ and $[\iota_2, \iota_4]$ belong to $\mathrm{IA}(W_{n,d})$. Hence we obtain the required result.

\section{The case where $n=2$}\label{S-n2}

In this section, we show that $\mathcal{K}_2$ is not finitely generated. More precisely, we show that
$\mathcal{K}_2$ is isomorphic to the second derived subgroup $[[F_2, F_2], [F_2, F_2]]$ of $F_2$.

\vspace{0.5em}

To begin with, we remark that $\mathrm{IA}_2$ is equal to the inner automorphism group $\mathrm{Inn}\,F_2$ of $F_2$ due to Nielsen \cite{Ni0}.
Furthermore, Bachmuth \cite{Ba1} showed that $\mathrm{IA}_2^M$ coincides with $\mathrm{Inn}\,F_2^M$. Hence we see
\[ \mathcal{K}_2 = \mathrm{Ker}(\mathrm{Inn}\,F_2 \xrightarrow{\mu_2} \mathrm{Inn}\,F_2^M). \]
In general, for any group $G$, a homomorphism from $G$ onto the inner automorphism group $\mathrm{Inn}\,G$ of $G$ defined by
\[ g \mapsto \iota_g = (x \mapsto g^{-1} x g, \hspace{1em} x \in G) \]
induces an isomorphism $\mathrm{Inn}\,G \cong G/Z(G)$ where $Z(G)$ denotes the center of $G$. Since the center of $F_2$ is trivial, we have
$\mathrm{Inn}\,F_2 \cong F_2$.

\vspace{0.5em}

Here we show that $Z(F_n^M)$ is trivial for any $n \geq 2$.
By a result of Magnus \cite{Mg0}, $F_n^M$ is faithfully represented in the group of $2 \times 2$ matrices of the form
\[ \begin{pmatrix} x & L \\ 0 & 1 \end{pmatrix} \]
where $x \in H$ and $L$ is a linear form $v_1 t_1 + \cdots + v_n t_n$ in the indeterminates $t_1, \ldots, t_n$ with coefficients in $\Z[H]$.
Here we identify $F_n^M$ with the image of it by the representation. The image of $x_i \in F_n^M$ by this representation is given by
\[ \begin{pmatrix} x_i & t_i \\ 0 & 1 \end{pmatrix}. \]
Let
\[ \begin{pmatrix} x & v_1 t_1 + \cdots + v_n t_n \\ 0 & 1 \end{pmatrix} \in Z(F_n^M). \]
Then for any $1 \leq i \leq n$, we have
{\small
\[\begin{split}
   \begin{pmatrix} x & v_1 t_1 + \cdots + v_n t_n \\ 0 & 1 \end{pmatrix} \begin{pmatrix} x_i & t_i \\ 0 & 1 \end{pmatrix}
    & = \begin{pmatrix} x x_i & v_1 t_1 + \cdots + (v_i+x) t_i + \cdots + v_n t_n \\ 0 & 1 \end{pmatrix}, \\
   \begin{pmatrix} x_i & t_i \\ 0 & 1 \end{pmatrix} \begin{pmatrix} x & v_1 t_1 + \cdots + v_n t_n \\ 0 & 1 \end{pmatrix}
    & = \begin{pmatrix} x_i x & (x_i v_1) t_1 + \cdots + (x_i v_i + 1) t_i + \cdots + (x_i v_n) t_n \\ 0 & 1 \end{pmatrix},
  \end{split}\]
}
and hence,
\[\begin{split}
    x_i v_j & = v_j, \hspace{1em} 1 \leq j \leq n, \,\,\, j \neq i, \\
    v_i + x & = x_i v_i +1.
  \end{split}\]
Since $\Z[H]$ is a domain, we see $v_j=0$ for $j \neq i$ from the former equation above. Since we can take $1 \leq i \leq n$ arbitrarily,
we obtain $v_1 = \cdots = v_n = 0 \in \Z[H]$. Then $x=1 \in H$. Therefore we see $Z(F_n^M)$ is trivial.

\vspace{0.5em}

By the argument above, we verify that $\mathrm{Inn}\,F_2^M \cong F_2^M$, and that $\mathcal{K}_2$ is the kernel
of the natural projection $F_2 \rightarrow F_2^M$. That is,
\begin{pro}
$\mathcal{K}_2 = [[F_2, F_2], [F_2, F_2]]$.
\end{pro}

\vspace{0.5em}

By an easy argument in the Combinatorial group theory, we see that $[[F_2, F_2], [F_2, F_2]]$ is a free group of infinite rank as follows.
First, we consider
\begin{lem}\label{L-F}
If $F$ is a free group of countably infinite rank, so is its commutator subgroup $[F,F]$. 
\end{lem}
\textit{Proof.}
Let $Y := \{ y_1, y_2, \ldots \}$ be the set of generators of $F$.
Then, applying the Reidemeister-Schreier method to $Y$ and a Schreier-transversal 
\[ S := \{ y_1^{e_1} y_2^{e_2} \cdots \,|\, e_i \in \Z, \,\,\, \text{all but finite $e_i$s are $0$.} \} \]
of $[F, F]$ of $F$,
we see that $[F, F]$ is a free group with basis
\[ \{ (s, y)=sy(\overline{sy})^{-1} \,|\, s \in S, \,\,\, y \in Y, \,\,\, (s, y) \neq 1 \}. \]
Hence, it suffices to show the set above contains infinitely many elements.

\vspace{0.5em}

Let us consider $y=y_1$ and $s=y_1^{e_1} y_2$ for any $\Z$. Then we have
\[ (s,y) = y_1^{e_1} y_2 y_1 y_2^{-1} y_1^{-(e_1+1)} \neq 1 \]
and
\[ y_1^{e_1} y_2 y_1 y_2^{-1} y_1^{-(e_1+1)} \neq y_1^{e_1'} y_2 y_1 y_2^{-1} y_1^{-(e_1'+1)} \]
if $e_1 \neq e_1'$. This completes the proof of Lemma {\rmfamily \ref{L-F}}. 

\vspace{0.5em}

In general, since the commutator subgroup $[F_n,F_n]$ of the free group $F_n$ is a free group of countably infinite rank for $n \geq 2$,
using Lemma {\rmfamily \ref{L-F}}, we conclude
\begin{pro}\label{P-main}
The group $\mathcal{K}_2$ is a free group of countably infinite rank.
\end{pro}

\section{Acknowledgments}\label{S-Ack}

This research is supported by JSPS Research Fellowship for Young Scientists and the Global COE program at Kyoto University.

\end{document}